\documentclass[11pt]{article}

\usepackage{amssymb,latexsym,amsmath}

\begin{document}

\newtheorem{thm}{Theorem}[section]

\newtheorem{prop}[thm]{Proposition}

\newtheorem{lema}[thm]{Lemma}

\newtheorem{cor}[thm]{Corollary}

\newtheorem{defi}[thm]{Definition}

\newtheorem{hypo}{Brittle fracture hypothesis}[section]

\newtheorem{ehypo}[thm]{Equilibrium hypothesis (EH)}
\newtheorem{sehypo}[thm]{Strong equilibrium hypothesis (SEH)}

\newtheorem{rk}[thm]{Remark}

\newtheorem{exempl}{Example}[section]

\newenvironment{exemplu}{\begin{exempl}  \em}{\hfill $\surd$

\end{exempl}}

\newcommand{\ea}{\mbox{{\bf a}}}
\newcommand{\eu}{\mbox{{\bf u}}}
\newcommand{\ep}{\mbox{{\bf p}}}
\newcommand{\ed}{\mbox{{\bf d}}}
\newcommand{\eD}{\mbox{{\bf D}}}
\newcommand{\eK}{\mathbb{K}}
\newcommand{\eL}{\mathbb{L}}
\newcommand{\eB}{\mathbb{B}}
\newcommand{\ueu}{\underline{\eu}}
\newcommand{\ueo}{\overline{u}}
\newcommand{\oeu}{\overline{\eu}}
\newcommand{\ew}{\mbox{{\bf w}}}
\newcommand{\ef}{\mbox{{\bf f}}}
\newcommand{\eF}{\mbox{{\bf F}}}
\newcommand{\eC}{\mbox{{\bf C}}}
\newcommand{\en}{\mbox{{\bf n}}}
\newcommand{\eT}{\mbox{{\bf T}}}
\newcommand{\eV}{\mbox{{\bf V}}}
\newcommand{\eU}{\mbox{{\bf U}}}
\newcommand{\ev}{\mbox{{\bf v}}}
\newcommand{\eve}{\mbox{{\bf e}}}
\newcommand{\uev}{\underline{\ev}}
\newcommand{\eY}{\mbox{{\bf Y}}}
\newcommand{\eP}{\mbox{{\bf P}}}
\newcommand{\eS}{\mbox{{\bf S}}}
\newcommand{\eJ}{\mbox{{\bf J}}}
\newcommand{\leb}{{\cal L}^{n}}
\newcommand{\eI}{{\cal I}}
\newcommand{\eE}{{\cal E}}
\newcommand{\hen}{{\cal H}^{n-1}}
\newcommand{\eBV}{\mbox{{\bf BV}}}
\newcommand{\eA}{\mbox{{\bf A}}}
\newcommand{\eSBV}{\mbox{{\bf SBV}}}
\newcommand{\eBD}{\mbox{{\bf BD}}}
\newcommand{\eSBD}{\mbox{{\bf SBD}}}
\newcommand{\ecs}{\mbox{{\bf X}}}
\newcommand{\eg}{\mbox{{\bf g}}}
\newcommand{\paromega}{\partial \Omega}
\newcommand{\gau}{\Gamma_{u}}
\newcommand{\gaf}{\Gamma_{f}}
\newcommand{\sig}{{\bf \sigma}}
\newcommand{\gac}{\Gamma_{\mbox{{\bf c}}}}
\newcommand{\deu}{\dot{\eu}}
\newcommand{\dueu}{\underline{\deu}}
\newcommand{\dev}{\dot{\ev}}
\newcommand{\duev}{\underline{\dev}}
\newcommand{\weak}{\rightharpoonup}
\newcommand{\weakdown}{\rightharpoondown}
\renewcommand{\contentsname}{ }

\title{Hamiltonian inclusions with convex dissipation with a view 
towards  applications}
\author{Marius Buliga\footnote{"Simion Stoilow" Institute of Mathematics of the Romanian Academy,
 PO BOX 1-764,014700 Bucharest, Romania, e-mail: Marius.Buliga@imar.ro }}
\date{This version: 08.06.09}


\maketitle

\begin{abstract}
We propose a generalization of Hamiltonian mechanics, as a Hamiltonian inclusion
with convex dissipation function. We obtain a dynamical version of the approach
of Mielke to quasistatic rate-independent processes.  
Then we show that a class of models of dynamical brittle damage can be
formulated in this setting. 
\end{abstract}

 {\bf MSC 2000:} 34G25; 70H05; 74R05

{\bf Keywords:}  Hamiltonian methods;  nonlocal damage; convex dissipation


\section{Introduction}

We are interested in the modification of the Hamiltonian formalism by adding 
the subdifferential of a convex dissipation function. In the Lagrangian formalism
this can be traced back to Rayleigh and Kelvin (cf. Thomson and Tait 
\cite{thomson} or Chetayev \cite{chetayev}). For the case of autonomous 
 Hamiltonian systems with a Rayleigh dissipation function added see the paper 
 of Bloch, Krishnaprasad, Marsden and Ratiu \cite{bloch}. Perturbation analysis of Hamiltonian
 systems is an old and much explored subject, which is beyond our scope in this 
 paper. 
 
 Closer to our interests is Mielke \cite{mielke} theory of quasistatic
 rate-independent processes. In fact one of our purposes is to reformulate 
 Mielke theory in a dynamical context. From this point of view a dissipation
 perturbed Hamiltonian approach seems the most economical. 
 
 From the viewpoint of multivalued analysis, many generalizations of Hamiltonian
 and Lagrangian mechanics have been considered, like for example Rockafellar 
\cite{rocka}, Aubin, Cellina and Nohel \cite{aubin} or Clarke \cite{clarke}. 
The problem of solving a subdifferential inclusion of the type 
(\ref{mainsymp}) for a 1-homogeneous dissipation function seems to be new. 
As a general problem this subdifferential inclusion seem to fall in the class 
of problems studied in the viabilty theory, \cite{aubin} or the more recent 
\cite{aubin2}, but the mathematical results 
in these papers do not apply here mainly because the dissipation is 
 1-homogeneous. 
 
More specifically, concerning the particular form -- (\ref{newqq}) coupled with purely
Hamiltonian equations (\ref{h1loceq}), (\ref{h2loc})-- of this subdifferential
inclusion, which is relevant for damage models in continuum mechanics, it seems
that  there are no mathematical results which  could be applied to 
this problem as a perturbed Hamiltonian problem. We thank to one of the 
anonymous referees for pointing us to the paper \cite{roubi}. From our viewpoint 
the results of this paper can be seen as leading to a solution of  our
problem, studied from the Lagrangian side, that is after reformulating it as a
generalized Euler-Lagrange equation. Nevertheless we think that the Hamiltonian
structure of this problem may lead to interesting discretization algorithms,
maybe based on symplectic integrators, which are known to handle correctly the 
energy balance even in the discretized form.

\paragraph{Outline of the paper.} In section 2 we propose and begin the study of
 a generalized 
Hamiltonian formalism, in the form of a subdifferential inclusion 
using a convex dissipation function. In
section 3 we show that Mielke's theory of quasistatic evolutionary processes is 
the quasistatic approximation of the formalism presented here. As an application, 
in section 4 we use the formalism for a energy of the Ambrosio-Tortorelli type and 
a 1-homogeneous dissipation function and we obtain a dynamical model of 
brittle damage which may be of interest in continuum media mechanics.

\paragraph{Acknowledgements.} This work has been done during two visits to LMT Cachan, due to the 
kind invitation of Olivier Allix. I want to thank him for  introducing me into the subject of delayed 
damage models, as well as for many constructive discussions.

\section{Generalized Hamiltonian equations with convex dissipation}
\label{proposed}

In the Lagrangian formalism we study the evolution of a system described by 
a variable $q$, which satisfies  the Euler-Lagrange equation associated to a 
function $\displaystyle L = L(t,q, \dot{q})$:  
\begin{equation}
 D_{q} L(t,q , \dot{q}) \, - \, \frac{d}{dt} D_{\dot{q}} L(t, q, \dot{q}) = 0
\quad \quad . 
\label{eleq}
\end{equation}
The function $L$ is called a Lagrangian and in many situation it has the form
\begin{equation}
L(t, q,  \dot{q}) \, = \, \hat{T}(\dot{q}) - \mathcal{E}(t,q)
\label{forml}
\end{equation}
 where:  $\displaystyle \hat{T}$ represents the kinetic energy,  is a smooth strictly convex smooth function 
 (for example quadratic, positive definite), and $\mathcal{E}$ is a potential energy 
 or stored energy. 
 
 In the formalism of Hamiltonian mechanics we double the variables: the system is 
 described by a pair $(q,p)$ where $p$ has the meaning of a momentum associated with
 $q$.  Instead of the Euler-Lagrange equation, the following   system of equations 
is used: 
\begin{equation}
 \left\{ 
 \begin{array}{rcl}
 - \dot{p} & \in & D_{q} H(t, q,p) \,  \\
 \dot{q} & = & D_{p} H(t, q,p)
 \end{array} \right. 
 \label{hhameq}
 \end{equation}
 The function $\displaystyle H = H(t,q, p)$ is called  a Hamiltonian.

 Consider for simplicity that 
$q, p \in \mathcal{H}$, where $\mathcal{H}$ is a Hilbert space with scalar product 
$\displaystyle \left( \cdot , \cdot \right)$. The equations of Hamiltonian mechanics
can be written in a compact form if we use the notations $\displaystyle 
z = (q,p) \in \mathcal{H} \times \mathcal{H}$, $J(z) = J(q,p) = (- p, q)$: 
\begin{equation}
J \, \dot{z} \, - \, D_{z} H(t,z) \, = \, 0 
\label{hhhameq}
\end{equation}
 
 In particular the Hamiltonian may take the form  
 \begin{equation}
 H(t,q,p) = T(p) \, + \, \mathcal{E}(t, q)
 \label{formh}
 \end{equation}
  where $T$ represents again the kinetic energy, this time expressed as a 
 function of $p$.

  In this case the two formalisms are equivalent if we take 
 $T$ to be the Fenchel conjugate of $\hat{T}$: 
 $$T(p) \, = \, \sup \left\{ \left( p,q \right) - \hat{T}(q) \mbox{ : } q \in
 \mathcal{H} \right\}$$

 \subsection{Introducing dissipation}
 
 Consider  a "dissipation function" $\mathcal{R}(q, \dot{q})$,  convex in the
 second argument, and a Lagrangian function which is a sum of kinetic and potential
 energies. In the particular case of 
 Rayleigh dissipation the  function $D$ has the form
 $$\mathcal{R}(q, \dot{q}) \, = \, \frac{1}{2} \| \dot{q} \|^{2}$$ 
 where $\| \cdot \|$ is a norm function. Then the Euler-Lagrange equations
 perturbed with the dissipation function $D$ are, by definition: 
 \begin{equation}
 D_{q} L(t,q , \dot{q}) \, - \, \frac{d}{dt} D_{\dot{q}} L(t, q, \dot{q}) \, \in \, 
 \partial_{\dot{q}} \mathcal{R}( q, \dot{q}) \quad \quad .
 \label{lag1}
 \end{equation}
 where the $\displaystyle \partial$ symbol denotes the subdifferential 
 from convex analysis. 
 
 The Hamiltonian side of (\ref{lag1}) is then 
 \begin{equation}
 \left\{ 
 \begin{array}{rcl}
 - \dot{p} & \in & D_{q} H(t, q,p) \, + \,  \partial_{\dot{q}} \mathcal{R}( q, \dot{q}) \\
 \dot{q} & = & D_{p} T(p) 
 \end{array} \right. 
 \label{ham1}
 \end{equation}
 
 This motivates us to  propose the following generalization of the Hamiltonian 
 equations (\ref{hhhameq}) 
in the form of a  subdifferential inclusion: 
\begin{equation}
J \, \dot{z} \, - \, D_{x} \, H(t,z) \, \in \, \partial_{\dot{z}} \,
\mathcal{R}(z,\dot{z})
\label{mainrel}
\end{equation} 
where $\displaystyle  \partial_{\dot{z}} \, \mathcal{R}(z,\dot{z})$ is the subdifferential 
of $D$ with respect to $\displaystyle \dot{z}$: 
\begin{equation}
\partial_{\dot{z}} \, \mathcal{R}(z,\dot{z}) \, = \, \left\{ (\bar{q}, \bar{p}) \in \mathcal{H} \times
\mathcal{H} \mbox{ : } \forall z' = (q',p') \in \mathcal{H} \times
\mathcal{H} \right. 
\label{defsubd}
\end{equation}
$$\left. \mathcal{R}(z, \dot{z} + z') \, \geq \, \mathcal{R}(z, \dot{z}) \, + \, \left( \bar{q}, q'
\right) \, + \, \left( \bar{p} , p' \right) \right\} \quad . $$

 We shall then be interested in the following particular case: suppose that we 
 have a decomposition of the state variable $\displaystyle q = (q_{1}, q_{2})$ 
 into a non-dissipative $\displaystyle q_{1}$ variable and a dissipative 
 $\displaystyle q_{2}$ variable. Then the momentum variable $p$ decomposes as 
 $\displaystyle p = (p_{1}, p_{2})$. The Hamiltonian function $H$ is taken as follows
 
 \begin{equation}
 H(t, q_{1}, q_{2}, p_{1}, p_{2}) = \mathcal{K}(p_{1}) + \frac{1}{2}
  \langle A p_{2} , p_{2} \rangle + \mathcal{E}(t, q_{1}, q_{2}) 
 \label{choiceh}
 \end{equation}
 where $\mathcal{K}$ is the kinetic energy energy associated to the variable 
 $\displaystyle p_{1}$, $A$ is a strictly positive definite 
  symmetric operator  and $\mathcal{E}$ is a stored energy function. 
 The dissipation function takes the form 
 \begin{equation}
 \mathcal{R}(q_{1}, q_{2}, \dot{q}_{1}, \dot{q}_{2}) = \rho(\dot{q}_{2})
 \label{choiced}
 \end{equation}
 with $\rho$ a convex function. 
 
 With these choices of functions $H$ and $D$ the system of equations (\ref{ham1}) becomes: 
 \begin{equation}
 \left\{ 
 \begin{array}{rcl}
 - \dot{p_{1}} & = & D_{q_{1}} \mathcal{E}(t, q_{1},q_{2}) \\
 \dot{q_{1}} & = & D_{p} \mathcal{K}(p_{1}) \\
 - \dot{p_{2}} & \in & D_{q_{2}} \mathcal{E}(t, q_{1},q_{2}) \, + \,  
 \partial \rho( \dot{q_{2}}) \\
 \dot{q_{2}} & = & A \, p_{2} \quad \quad .
 \end{array} \right. 
 \label{ham2}
 \end{equation}
We can see the first two equations as a Hamiltonian evolution of 
the variables $\displaystyle (q_{1}, p_{1})$ which has $\displaystyle 
(q_{2}, p_{2})$ as control parameters, coupled with a pair of evolution
equations (the last two equations in (\ref{ham2})) for the control parameters. 
These last two equations can be see as a differential inclusion:  
  \begin{equation}
  - \dot{p_{2}} \, - \, D_{q_{2}} \mathcal{E}(t, q_{1},q_{2}) \, \in \,  
 \partial \rho( A^{-1} p_{2}) \quad \quad .
 \end{equation}

 Interesting particular cases of dissipation function $\rho$ are: 
\begin{enumerate}
\item[(a)] $\rho =0$, no dissipation, this corresponds to classical Hamiltonian
equations, 
\item[(b)] $\displaystyle \rho(\dot{q}) \, = \, \frac{1}{2} \left( 
\dot{q}, \dot{q}\right)$, (where $\left( \cdot, \cdot \right)$ is a scalar
product),  which can be traced back to the  Rayleigh dissipation function, 
 \item[(c)] $\displaystyle \rho(\dot{q}) \, = \,  \| 
\dot{q} \|$, where $\| \cdot \|$ is a Banach space norm, or  a more general 
1-homogeneous convex function which, as we shall explain, 
is related to the approach of Mielke and collaborators -- Mielke, Theil 
\cite{mielketh99},   Mielke, Theil and 
Levitas \cite{mielkethl},  \cite{mielke} --  to quasistatic rate-independent 
evolutionary processes.
\end{enumerate}

\subsection{The formalism in topological vector spaces}

We shall precisely formulate relation (\ref{mainrel}) for a pair of locally
compact topological vector spaces in duality. in particular this will cover 
the cases of Banach or Hilbert spaces.

$X$ and $Y$ are topological, locally convex, real vector spaces of dual 
variables $x \in X$ and $y \in Y$, with the duality product 
$\langle \cdot , \cdot \rangle : X \times Y \rightarrow \mathbb{R}$. 
We shall suppose that $X, Y$ have topologies compatible with the duality 
product, that is: any  continuous linear functional on $X$ (resp. $Y$) 
has the form $x \mapsto \langle x,y\rangle$, for some $y \in Y$ (resp. 
$y \mapsto \langle x,y\rangle$, for some  $x \in X$).

In this frame we don't have scalar products, neither an equivalent of the 
linear transformation $J$, therefore we start by introducing natural notations 
which make sense in this generality. 

We want to study generalized Hamiltonian evolutions in the space $X \times Y$. 
For a general element of $X \times Y$ we shall use the notation $z = (x,y)$, or
similar.  

In order to properly formulate Hamiltonian equations or inclusions we need: 
a symplectic form, a Poisson bracket and a notion of subdifferential adapted in
this setting. These will be the most natural objects one may think about and
they were used many times before.

We shall use notations familiar in symplectic geometry, namely: $\omega$ for the
symplectic form, $\left\{ \cdot , \cdot \right\}$ for the Poisson bracket, 
$\displaystyle X_f$ for the symplectic gradient of  the function 
$f: X \times Y \rightarrow \mathbb{R}$ (if the linear $J$ is available then 
$\displaystyle X_f \, = \, - J \, D f$, where $D f$ is the differential of $f$). Instead of
the usual  subdifferential of a convex function $F$ we shall use a "symplectic 
subdifferential" $X \, F$. In the usual setting in Hilbert spaces we have 
$J \, X \, F \, = \, \partial \, F$, where $\partial F$ is the well known 
subdifferential from convex analysis. In this general setting the definition 
of $X \, F$ is obtained from the definition of $\partial F$ by replacing scalar
products with the symplectic form.

Remark however that in this general setting the symplectic form and Poisson
bracket have to be understood in a weaker sense than usual, let's say on a
finite dimensional symplectic manifold. Indeed, a symplectic
form is a non-degenerated 2-form which is closed (we renounce to the condition of being closed);
a Poisson bracket is a Lie bracket over a algebra of functions, with supplementary
properties,  while here the "Poisson bracket" we define sends a pair of 
differentiable functions from $Der(X,Y)$ to a function which is not
differentiable a priori. 

Let us proceed with the introduction of the necessary objects.

\begin{defi}

The  space $X \times Y$  is endowed with a  symplectic form: 
for any $z'=(x',y')$, 
$z"=(x",y")$ we define the bilinear and anti-symmetric form  
$$\omega ( z', z" ) \, = \langle x' , y" \rangle - \langle x" , y' \rangle \quad
\quad . $$

$\displaystyle Der(X,Y)$ is the linear space of  functions 
$f: X\times Y \rightarrow \mathbb{R}$ which are continuously 
differentiable in each argument
in the following sense: there are continuous functions $\displaystyle 
D_{x} f : X\times Y \rightarrow Y$ and $\displaystyle 
D_{y} f : X\times Y \rightarrow X$ such that 
for any $(x,y) \in X \times Y$ and 
\begin{enumerate}
\item[(a)]  for all 
$y' \in Y$ we have 
$$\lim_{\varepsilon \rightarrow 0} \frac{1}{\varepsilon} \left[ f(x, y +
\varepsilon y') - f(x,y) \right] \, = \, \langle D_{y} f (x,y) , y' \rangle $$
\item[(b)]  for all $x' \in X$ we have 
$$\lim_{\varepsilon \rightarrow 0} \frac{1}{\varepsilon} \left[ f(x +
\varepsilon x', y ) - f(x,y) \right] \, = \, \langle x', D_{x} f (x,y)\rangle $$
\end{enumerate}
The symplectic gradient of $f \in Der(X,Y)$ is the function 
$\displaystyle X_f : X \times Y \rightarrow X \times Y$ defined by 
$$X_f (x,y) \, = \, (D_{y} f (x,y) , - D_{x} f(x,y) ) \quad \quad . $$

The Poisson bracket is the bilinear, antisymmetric form 
$$\displaystyle \left\{ \cdot , \cdot
\right\} : Der(X,Y) \times Der(X,Y) \rightarrow \mathbb{R}^{X\times Y}$$ 
defined by: 
$\displaystyle \left\{ f , g \right\} \, = \, \omega \left( X_f , X_g \right)$. 
\label{defobj}
\end{defi}

\begin{defi}
Let $F: X \times Y \rightarrow \mathbb{R}$ be a convex lsc function. The
symplectic subdifferential of $F$ is the multivalued function which sends 
$z = (x,y) \in X \times Y$ to the set 
$$\displaystyle X \, F (z) \, = \, \left\{ 
z' \in X \times Y \mbox{ : } \forall \, z" \in X \times Y \quad  
F(z+z") \,  \geq \,  F(z) \, + \, 
\omega (z' , z")\right\} $$ 
\label{defssub}
\end{defi}

Remark that if $F \in Der(X,Y)$ and convex then we have 
$\displaystyle X \, F \, = \, \left\{ X_{F} \right\}$. Indeed, if we use 
$\displaystyle z' = X_{F} = (D_{y} F(x,y) , - D_{x} F(x,y))$ in the 
definition \ref{defssub} of the symplectic differential we get 
$$F(z + z") \, \geq \,   F(z)  \, + \, \langle D_{y} F(x,y) , y" \rangle \, + \, 
\langle x" ,  D_{x} F(x,y) \rangle $$ 
which is true due to the convexity of $F$. Therefore $\displaystyle X_{F} (x,y) \in 
X \, F (x,y)$. The converse implication, that is $\displaystyle z' \in X \, F (x,y)$ 
implies $\displaystyle z' = X_{F}(x,y)$,  is true by standard arguments of 
convex analysis. 

We propose the following generalization of Hamiltonian evolution. 

\begin{defi}
Let $H : [0,T] \times X \times Y \rightarrow \mathbb{R}$ such that for all 
$t \in [0,T]$ we have $H(t, \cdot) \in Der(X,Y)$, and $\displaystyle 
D: (X \times Y)^{2} \rightarrow \mathbb{R} \cup \left\{ + \infty \right\}$ be a
"dissipation function" with the properties: 
\begin{enumerate}
\item[(a)] for any $z', z" \in X \times Y$  we have $\mathcal{R}(z', z") \geq 0$ and 
$\mathcal{R}(z' ,0) = 0$, 
\item[(b)] for any $z \in X \times Y$ the function $\mathcal{R}(z,  \cdot)$ is convex, 
lsc. 
\end{enumerate}
Then a curve $z: [0,T] \rightarrow X \times Y$ is a solution of the 
evolution problem with Hamiltonian
$H$ and dissipation $D$ if it is derivable for all $t \in [0, T]$ (with
differential denoted by $\displaystyle \dot{z}$) and it
satisfies the subdifferential inclusion: 
\begin{equation}
\dot{z}(t) \, - \, X_{H(t,\cdot)} (z(t)) \, \in \, X\, \left( \mathcal{R}(z(t), \cdot )
\right) (\dot{z}(t)) \quad \quad .  
\label{mainsymp}
\end{equation}
\label{mainevol}
\end{defi}

We can give an equivalent characterization for a solution, which later will lead
to a notion of weak solution. For any $f \in Der(X,Y)$ and any derivable curve 
$z: [0,T] \rightarrow X \times Y$ we denote by $f \circ z : [0,T] \rightarrow 
\mathbb{R}$ the function composition of $f$ and $z$, and by $\displaystyle 
\frac{d}{d \, t} \left[ f \circ z \right](t)$ the differential of this
composition. 

\begin{prop}
With the notations from definition \ref{mainevol}, $z$ is a solution of the
evolution problem if and only if for any $f \in Der(X,Y)$ and for any $t \in
[0,T]$  we have: 
\begin{equation}
\mathcal{R}(z(t), \dot{z}(t) - X_{f}(z(t))) \, \geq \, \mathcal{R}(z(t), \dot{z}(t)) \, + \, 
\frac{d}{d \, t} \left[ f \circ z \right](t) \, - \, \left\{ f , H(t, \cdot )
\right\} (z(t)) \quad \quad . 
\label{poissonfor}
\end{equation}
\label{poissonprop}
\end{prop}

\paragraph{Proof.} For any $f \in Der(X,Y)$ and any derivable curve 
$z: [0,T] \rightarrow X \times Y$ we have, by direct computation: 
\begin{equation}
\frac{d}{d \, t} \left[ f \circ z \right](t) \, - \, \left\{ f , H(t, \cdot )
\right\} (z(t)) \, = \, - \, \omega \left( \dot{z}(t) - X_{H(t,\cdot)}(z(t)) ,
X_{f}(z(t)) \right) \quad \quad . 
\label{useit}
\end{equation}
Let $z$ be a solution of the evolution problem. We choose then in (\ref{mainsymp}) 
$\displaystyle z" = - X_{f}(z(t))$ and use (\ref{useit}) to get
(\ref{poissonfor}). 

Conversely, suppose that the curve $z$ satisfies (\ref{poissonfor}). For any 
$\displaystyle z" \in X \times Y$ let us define $f \in Der(X,Y)$ by 
$f(z) = \omega(z, z")$. It is easy to see then that $\displaystyle 
X_{f} = - z"$,  that $\displaystyle \frac{d}{d \, t} \left[ f \circ z \right](t)
= \omega (\dot{z}(t), z")$ and that $\displaystyle \left\{ f , H(t, \cdot )
\right\} (z(t)) \, = \, \omega \left( X_{H(t, \cdot)}, z"\right)$. In conclusion
the relation (\ref{poissonfor}) for this choice of the function $f$ becomes 
the relation (\ref{mainsymp}) for $z"$. \quad $\square$

\vspace{.5cm}

It is visible that the functions $f \in Der(X,Y)$ play the role of test
functions in (\ref{poissonfor}). 
Let us consider  curves 
$f: [0,T] \rightarrow Der(X,Y)$, which are smooth in the sense that for any 
$t \in [0,T]$ there exists $\displaystyle \frac{\partial}{\partial \, t} 
f(t,z)$. We  suppose that the Hamiltonian $H: [0,T]
\rightarrow Der(X,Y)$ is such a curve. For an arbitrary $t \in [0,T]$, 
at each $\tau \in [0,t]$ we put 
$f(\tau,\cdot)$ in the relation (\ref{poissonfor}) and then integrate with respect 
to $\tau \in [0,t]$. We  obtain the following
relation: 
$$\int_{0}^{t} \mathcal{R}(z(\tau), \dot{z}(\tau) - X_{f(\tau,\cdot)}(z(\tau))) \mbox{
d}\tau \, 
\geq \, \int_{0}^{t} \mathcal{R}(z(\tau), \dot{z}(\tau)) \mbox{ d}\tau \, + \, \quad \quad \quad \quad \quad
$$ 
\begin{equation}
\quad \quad \quad \quad \quad \quad \quad \quad \quad + \, f(t, z(t)) \, - \, f(0, z(0)) \, - \,
\label{weakform}
\end{equation}
$$ \quad \quad \quad \quad \quad \quad  - \, \int_{0}^{t} \left[ \frac{\partial}{\partial \, t} 
f(\tau,z(\tau)) \, + \, \left\{ f(\tau, \cdot) , H(\tau, \cdot)
\right\}(z(\tau)) \right] \mbox{ d}\tau \quad . $$
The relation (\ref{weakform}) makes sense  if $z$ is differentiable almost
everywhere  and 
\begin{equation}
 \int_{0}^{T} \mathcal{R}(z(\tau),
\dot{z}(\tau)) \mbox{ d}\tau \, < \, + \infty
\label{cond1}
\end{equation}
\begin{equation}
\int_{0}^{t} \left[ \frac{\partial}{\partial \, t} 
f(\tau,z(\tau)) \, + \, \left\{ f(\tau, \cdot) , H(\tau, \cdot)
\right\}(z(\tau)) \right] \mbox{ d}\tau \, < \, + \infty \quad . 
\label{cond2}
\end{equation}
This is leading us to the following definition of weak solution.

\begin{defi}
Let $\mathcal{A}$ be a given vector space of smooth curves 
$f: [0,T] \rightarrow Der(X,Y)$ such that the Hamiltonian $H: [0,T]
\rightarrow Der(X,Y)$  belongs to $\mathcal{A}$. Then let 
$\mathcal{S}(D, \mathcal{A})$ be the space of all curves 
$\displaystyle z: [0,T] \rightarrow X \times Y$ which are almost everywhere 
derivable, such that $Diss(z, [0,T]) < + \infty$ and such that (\ref{cond2}) is true for any $f \in \mathcal{A}$. 

A curve $z \in \mathcal{S}(D, \mathcal{A})$ is a weak solution of the evolution
problem if for almost any $t \in [0,T]$ the inclusion (\ref{mainsymp}) is true. 

Let $z \in \mathcal{S}(D, \mathcal{A})$ be a weak solution. The dissipation 
along  this solution is by definition the function: 
\begin{equation}
\eta(t) \, = \, \, \int_{0}^{t} \omega \left( \dot{z}(\tau) ,X_{H(\tau, \cdot)} 
(z(\tau)) \right) \mbox{ d}\tau  \quad \quad .
\label{dishameq}
\end{equation}

\label{defweaksol}
\end{defi}

\begin{prop}
Let $z \in \mathcal{S}(D, \mathcal{A})$ be a weak solution
and $\eta$ the associated dissipation. Then for any $t \in [0,T]$ we have 
$$\eta(t) \, \geq \, \int_{0}^{t} \mathcal{R}(z(\tau), \dot{z}(\tau)) \mbox{
d}\tau \quad . $$ 
\label{pdis}
\end{prop}

\paragraph{Proof.}
We shall use the inclusion (\ref{mainsymp}), which means that for any 
$\displaystyle z' \in X \times Y$ and for almost any $t \in [0,T]$ we have 
$$\mathcal{R}(z(t), \dot{z}(t) + z') \, \geq \, \mathcal{R}(z(t), \dot{z}(t)) \,
+ \, \omega ( \dot{z}(t) - X_{H(t, \cdot)}(z(t)) , z') \quad . $$
If we take   for almost any $\tau \in [0,t]$ $\displaystyle z' = - \dot{z}(\tau)$ and use $\mathcal{R}(z, 0) = 0$
then we get 
$$\omega \left( \dot{z}(\tau) ,X_{H(\tau, \cdot)} 
(z(\tau)) \right) \, \geq \, \mathcal{R}(z(\tau), \dot{z}(\tau)) \, \geq \, 0
\quad . $$
The desired relation is obtained by integration.  \quad $\square$

\subsection{The 1-homogeneous case}

Suppose that $X$ is a Banach space and $Y = X^{*}$. Then $X \times Y$ is a
Banach space and the natural norm on $X \times Y$ induces a distance $d(z' , z")
= \| z' - z" \| $. 

Suppose moreover that for any $z \in X \times Y$ the dissipation function 
$D$ has the property that $\mathcal{R}(z, \cdot)$ is positively one-homogeneous. 
Then the 
dissipation function can be seen as a dissipation metric in the sense that 
it induces: 
\begin{enumerate}
\item[(a)] a "dissipation length" defined for any curve $z: [0,T] \rightarrow 
X \times Y$ which is almost everywhere differentiable by: 
$$L(z) \, = \, \int_{0}^{T} \mathcal{R}(z(t), \dot{z}(t)) \mbox{ d}t$$
The space of
curves with finite dissipation length is denoted with $\displaystyle 
W^{1,1}_{\mathcal{R}}(X \times Y)$. 
\item[(b)] a "dissipation distance" $\displaystyle D : (X \times Y)^{2}
\rightarrow \mathbb{R} \cup \left\{ + \infty \right\}$, where 
$D (z', z")$ is  defined as the infimum of
the dissipation lengths of all curves joining $z'$ and $z"$. 
\item[(c)] a "dissipation variation" defined for any curve $z: [0,T] \rightarrow 
X \times Y$ as: 
$$Diss(z, [0,T]) \, = \, \sup \left\{ \sum_{1}^{N} D (z(s_{j-1}), z(s_{j}))
\, \mid \, \mbox{all partitions of } [0,t] \right\} \quad . $$
$\displaystyle  BV_{\mathcal{R}}(X \times Y)$ denotes the
 space of curves with bounded dissipation variation. 
\end{enumerate}

The dissipation distance $\displaystyle D $ is not really a distance, because
it is not symmetric and it may take the value $+ \infty$. It satisfies
nevertheless the triangle inequality. 

The dissipation length and dissipation variation are defined in principle 
for different classes of curves, but in particular cases they are the same. 
All in all this is a  generalization of well-known facts in the analysis 
in metric spaces, see for the relevant results Gromov  chapter 3 \cite{gromov},
or Ambrosio, Gigli, Savar\'e  chapter 1\cite{amb}, which has been developed by 
Mielke and collaborators in the theory of rate-independent evolution systems 
(see section \ref{secmielke} for further details and references). Enough is to
mention that if $z$ is a curve which is differentiable almost everywhere and 
of finite dissipation length then its dissipation length equals the dissipation
variation.

In particular  then 
 any weak solution satisfies (\ref{weakform}) with the term 
$\displaystyle \int_{0}^{t} \mathcal{R}(z(\tau),
\dot{z}(\tau)) \mbox{ d}\tau $ replaced by $Diss(z, [0,t])$. If the class
$\mathcal{A}$ is sufficiently rich then satisfaction of (\ref{weakform}) will 
imply that $z$ is a weak solution. 

\begin{thm}  
If $\mathcal{R}(z', z") \, = \, \mathcal{R}(z', x")$ for any 
$z', z" \in X \times Y$ then for any weak solution 
$z: [0,T] \rightarrow X \times Y$ and for any $t \in 
[0,T]$ we have: 
\begin{equation}
H(0, z(0)) \, + \, \int_{0}^{t} \frac{\partial}{\partial \, t} H(\tau, z(\tau)) 
\mbox{ d} \tau \, = \, H(t, z(t)) \, + \, Diss(z, [0,t]) 
\label{enerbal}
\end{equation}
\label{tenerbal}
\end{thm}

\paragraph{Proof.} 
In relation (\ref{weakform}) let us take $f = \lambda H$ for an arbitrary 
$\lambda \in (-\infty, 1)$: 
$$\int_{0}^{t} \mathcal{R}(z(\tau), \dot{z}(\tau) - \lambda X_{H(\tau,\cdot)}(z(\tau))) \mbox{
d}\tau \, 
\geq \, Diss(z, [0,t]) \, + \, \lambda H(t, z(t)) \, - \, \lambda H(0, z(0)) \,
- \, 
$$ 
\begin{equation}
\quad \quad \quad \quad \quad \quad \quad \quad \quad  - \, \lambda 
\int_{0}^{t} \left[ \frac{\partial}{\partial \, t} 
H(\tau,z(\tau))  \right] \mbox{ d}\tau \quad .
\label{weakform1}
\end{equation}
In the hypothesis of the theorem if $z$ is a weak solution then it satisfies the
following: for almost any $t \in [0,T]$ and for 
 any $z" = (x", y") \in X \times Y$ 
$$\mathcal{R}(z(t),\dot{x}(t) + x") \, \geq \,  \mathcal{R}(z(t),\dot{x}(t)) \, + \, 
\langle \dot{x}(t) -  D_{y}H(t, \cdot) (x(t), y(t)), y" \rangle \, - \, $$
$$- \, 
\langle x", \dot{y}(t) + D_{x}H)t, \cdot) (x(t), y(t)) \rangle$$ 
It follows that for almost any $t \in [0,T]$ we have $\displaystyle 
\dot{x}(t) =  D_{y}H(t, \cdot) (x(t), y(t))$, therefore for almost any 
$\tau \in [0,t]$ we have: 
$$\mathcal{R}(z(\tau), \dot{z}(\tau) - \lambda X_{H(\tau,\cdot)}(z(\tau))) \, = \, 
\mathcal{R}(z(\tau), \dot{x}(\tau) - \lambda D_{y}H(\tau, \cdot) (x(\tau), y(\tau))) \, =
$$ 
$$ \, = \, \mathcal{R}\left(z(\tau), (1-\lambda) \left( \dot{x}(\tau) \right) \right) \,
= \, (1 - \lambda) \mathcal{R} (z(\tau), \dot{x}(\tau))$$ 
We return to (\ref{weakform1}), we use the information that we gained and the
equality between dissipation variation and dissipation distance and we obtain: 
for any $\lambda \in (-\infty, 1)$ we have: 
$$0 \, \geq \, \lambda \left[  Diss(z, [0,t]) \, + \, H(t, z(t)) \, - \, 
 H(0, z(0)) \, - \, \int_{0}^{t} \left[ \frac{\partial}{\partial \, t} 
H(\tau,z(\tau))  \right] \mbox{ d}\tau \right] \quad .$$
The arbitrary $\lambda$ can have any sign, therefore we deduce the desired
equality (\ref{enerbal}) from the previous inequality. \quad $\square$

\vspace{.5cm}

This theorem shows a great advantage of Hamiltonian formulations upon Lagrangian
formulations: a weak
Hamiltonian formulation naturally conserves quantities of interest, like the 
energy, while in Lagrangian formulations this has to be imposed by hand (which 
then leads to different  weak and energetic formulation). This can be stated in 
few words as: weak solutions of the Hamiltonian formulation are energetic 
solutions in the Lagrangian formulation.

\section{Connection with Mielke's theory of  quasistatic evolutionary processes}
\label{secmielke}

Consider a physical system with the state space  $\mathcal{Q}$. This space 
may have a manifold structure, or it may be a space  of functions $q: \Omega \rightarrow \mathcal{M}$,  with 
given regularity, where $\mathcal{M}$ is a manifold. In this case   
the bounded Lipschitz domain $\Omega$ represents the reference configuration  
of a continuous body. We shall  denote a generic 
point of $\mathcal{Q}$ by the letter $q$ and $\dot{q}$ denotes a vector in the tangent 
space to $\mathcal{Q}$ at $q \in Q$.  

For the first time in the proceedings paper Mielke, Theil 
\cite{mielketh99}, then in   Mielke, Theil and 
Levitas \cite{mielkethl},  the notion of a energetic solution  of a 
quasistatic evolutionary process was introduced,  
based on a energy function  
$$ \mathcal{E} : [0,T] \times \mathcal{Q} \rightarrow \mathbb{R} \cup \left\{ + \infty \right\} \quad , \quad \mathcal{E} = \mathcal{E}(t, q) $$
and a "dissipation metric" 
$$\mathcal{R}  : T \mathcal{Q} \rightarrow [0,+\infty] \quad , \quad \mathcal{R} = \mathcal{R}(q, \dot{q})$$
Here $\displaystyle T \mathcal{Q}  \, = \, \left\{Ê(q, \dot{q}) \, \mid \,  \dot{q} \in T_{q} \mathcal{Q} \right\}$ 
is the tangent space space to $\mathcal{Q}$ at $q \in Q$, in a generalized sense. 

The dissipation metric  is convex and lower semicontinuous with 
respect to  the second variable. For the case of rate-independent 
processes the dissipation metric  is  1-homogeneous (i.e. it can really be
interpreted as a metric).  
The {\bf force balance} equation is: 
\begin{equation}
0 \, \in \, \partial_{\dot{q}} \, \mathcal{R}(q, \dot{q}) \, + \, D_{q} \, \mathcal{E}(t,q)
\label{fmielke}
\end{equation}

To the dissipation metric $\mathcal{R}$ is associated a non symmetric dissipation distance 
$$D : \mathcal{Q} \times \mathcal{Q} \rightarrow [0, +\infty]$$
$$D(q_{1}, q_{2}) \, = \, \inf \left\{ \int_{0}^{1} \mathcal{R}(q(s), \dot{q}(s)) \mbox{Êd}s \, \mid \, 
q \in W^{1,1}([0,1], \mathcal{Q}) \, q(0) = q_{1}, \, q(1) = q_{2} \right\}$$

\begin{defi}
A evolution $q: [0,T] \rightarrow \mathcal{Q}$ is an energetic solution associated with $\mathcal{E}$ and 
$D$ if 
\begin{enumerate}
\item[(a)] the function $\displaystyle t \in [0,T] \mapsto \partial_{t} \mathcal{E}(t, q(t))$ belongs to 
$\displaystyle L^{1}((0,T))$, and for every $t \in [0,T]$ we have $\mathcal{E}(t,q(t)) < + \infty$, 
\item[(b)] the stability condition holds: for any $\hat{q} \in \mathcal{Q}$ 
$$ \mathcal{E}(t, q(t)) \, \leq  \, \mathcal{E}(t, \hat{q}) \, + \, D(q(t), \hat{q})$$
\item[(c)] the energy balance holds: 
$$ \mathcal{E}(t, q(t)) \, + \, Diss \, (q, [0,t]) \, = \, \mathcal{E}(0, q(0)) \, + \, \int_{0}^{T} \partial_{t} \mathcal{E}(s, q(s)) \mbox{ d}s$$
\end{enumerate}
where 
$$Diss  \, (q, [0,t]) \, = \, \sup \left\{ \sum_{1}^{N} D(q(s_{j-1}), q(s_{j})) \, \mid \, \mbox{Êall partitions of } [0,t] \right\}$$
\label{enmielke}
\end{defi}

We can recover the force balance equation (\ref{fmielke}) from the generalized 
Hamiltonian formalism with dissipation proposed in section \ref{proposed}.
Indeed, suppose that the state space of the physical system is 
$\mathcal{Q} = B$, a reflexive  Banach space. Consider the phase
space $\displaystyle X \times Y  =  B \times B^{*}$. A generic element of 
$z \in \mathcal{B}$ has the form $z = (q,p)$ with $q \in B$, $p \in B^{*}$. 

We shall take Hamiltonian and dissipation functions almost as in
(\ref{choiceh}), (\ref{choiced}). 
 The  Hamiltonian function $H$ has the form  $\displaystyle 
 H(t, q, p) = \mathcal{K}(p)  + \mathcal{E}(t, q)$ 
 where $\mathcal{K}$ is a smooth function (kinetic energy) and $\mathcal{E}$ is the  
 energy function of Mielke. We take a dissipation function  
$\displaystyle  \mathcal{R}(q,p, \dot{q}, \dot{p}) = \mathcal{R}(q,\dot{q})$ 
 with $\mathcal{R}$ the dissipation metric. 
 
 With these choices of functions $H$ and $D$ the equation (\ref{mainsymp}) takes the
 form: 
 \begin{equation}
 \left\{ 
 \begin{array}{rcl}
 - \dot{p} & \in & D_{q} \mathcal{E}(t, q) \, + \,  
 \partial_{\dot{q}} \mathcal{R}(q, \dot{q}) \\
 \dot{q} & = & D_{p} \mathcal{K}(p)  \quad \quad .
 \end{array} \right. 
 \label{ham22}
 \end{equation}
The quasistatic version of (\ref{ham22}) is just the force balance equation of
Mielke (\ref{fmielke}). We are also in the hypothesis of theorem \ref{tenerbal}.
If we neglect the inertial terms in (\ref{enerbal}) we obtain 
 the energy balance condition (c) from the 
definition of energetic solution \ref{enmielke}.

 Let us see what is the expression of the dissipation  
along a solution of (\ref{ham22}), as defined by (\ref{dishameq}). We have
$$\dot{\eta}(t) \, = \, - \, \langle D_{p} \mathcal{K}(p(t)), 
\dot{p}(t)\rangle \, - \, \langle \dot{q}(t), D_{q} \mathcal{E}(t,q(t))\rangle $$
 As in the proof of proposition \ref{pdis}, we arrive to the inequality 
$$0 \geq \, \mathcal{R}(q(t),\dot{q}(t)) \, + \, \langle D_{p} \mathcal{K}(p(t)), 
\dot{p}(t)\rangle \, + \, \langle \dot{q}(t), D_{q} \mathcal{E}(t,q(t))\rangle$$
therefore we get $\displaystyle \mathcal{R}(q(t),\dot{q}(t)) \leq
\dot{\eta}(t)$. We integrate this inequality and we obtain: 
$$\eta(t) \, \geq \, \int_{0}^{t} \mathcal{R}(q(s), \dot{q}(s)) \mbox{d}s$$
We finally obtain that $\displaystyle \eta(t) \, \geq \, 
D (q(0), q(t)) \, \geq \, 0$, 
which means that the dissipation $\eta$ along a solution of (\ref{ham22}) is
always greater or equal to the dissipation distance (in fact greater than the
dissipation length). 

\section{Application: a dynamical model of brittle damage using the
Ambrosio-Tortorelli functional}

Mielke and Roub\'{\i}\v{c}ek \cite{MR06b} proposed a rate-independent brittle damage model based 
on the theory of rate-independent evolutionary processes \cite{mielke}. 
The model of Mielke and Roub\'{\i}\v{c}ek is a quasistatic  particular case of the more general
 dynamical model of Stumpf and Hackl \cite{sthk}.  
 
By using the generalized Hamiltonian formalism we are able to obtain a dynamical
model of brittle damage, which is also a particular case of the general
dynamical model of Stumpf and Hackl. 

The model is based on a energy of Ambrosio-Tortorelli type and a dissipation 
function as in the model of Mielke and Roub\'{\i}\v{c}ek.

\subsection{The Ambrosio-Tortorelli functional}

Let $n\in \mathbb{N}$ be a strictly positive natural number and $\displaystyle \Omega \subset \mathbb{R}^{n}$ a 
bounded, open set, with piecewise smooth boundary.  
The Mumford-Shah functional \cite{MS} is 
\begin{equation}
E(u,S) \ = \  \int_{\Omega}  \frac{1}{2} \, K \, \mid \nabla u \mid^{2} \,    \, + \, \gamma \, \mathcal{H}^{n-1}(S) 
  \label{ms2}
\end{equation}
defined over all pairs $(u,S)$ such that $\displaystyle u \in \mathcal{C}^{1}( \Omega \setminus S ,  \mathbb{R})$. The set $S$ is a $n-1$-dimensional surface in $\displaystyle \mathbb{R}^{n}$, or a countable union of such surfaces.   In the case $n = 2$ this functional can be seen as the energy of a brittle body suffering a antiplane 
displacement $u$ and presenting a crack $S$.  

For $n=3$ the state of a brittle body is  described by a pair 
displacement-crack. $(\eu,S)$ is such a pair if $S$ is a crack (a 2D surface)  which appears in the body and  
$\displaystyle \eu \in \mathcal{C}^{1}( \Omega \setminus S ,  \mathbb{R}^{3})$ is a displacement of the broken body, that is  $\eu$ is smooth  in the exterior of the surface $S$, but it may have jumps over $S$.  
The total energy of a brittle  body is a Mumford-Shah functional of the form: 
\begin{equation}
E(\eu,S) \ = \ \int_{\Omega} w(\nabla \eu) \mbox{ d} x \, + \, G \, \mathcal{H}^{2}(S) \quad \quad .
  \label{ms3}
\end{equation}
The first term of the functional $E$ represents the elastic energy of the body with the
displacement $\eu$. The second term represents the energy consumed to produce the crack 
$S$ in the body. Here his energy is taken to be proportional with the area of the crack $S$ (technically 
this is the  2 dimensional Hausdorff measure of $S$), with the proportionality factor $G$, which is the Griffith 
constant. 

Starting with the foundational papers of Mumford, Shah \cite{MS},  De Giorgi, Ambrosio \cite{DGA}, Ambrosio \cite{A1}, \cite{A2},  
the development of  models of quasistatic brittle fracture  based on Mumford-Shah functionals  continues with Francfort, Marigo \cite{FMa}, 
\cite{FMa2}, Mielke \cite{mielke}, Dal Maso, Francfort, Toader, \cite{dft05}, 
Buliga \cite{bu}, \cite{bu3}, \cite{bu2}.  

All these models are based on a technique of time discretization followed by a sequence of incremental minimization problems. These models are either seen as applications of De Giorgi method of 
energy minimizing movements, or in the frame of the theory of Mielke of rate-independent evolutionary 
processes \cite{mielke}.

The functional 
\begin{equation}
E_{c}(u,d) \ = \  \int_{\Omega} \left\{ \phi(d) \, \frac{1}{2} \, K \, \mid \nabla u \mid^{2} \, + \, \frac{1}{2} \gamma \, c \, 
  \mid \nabla d \mid^{2} \, + \, \frac{\gamma}{2 c} \, d^{2} \right\} \,  
  \label{atf}
\end{equation}
was introduced by Ambrosio and Tortorelli \cite{ambtor}, as a variational 
approximation of the Mumford-Shah functional (\ref{ms2}).  Here $d$ is a field which approximates the 
characteristic function of a crack, that is $d: \Omega \rightarrow [0,1]$ 
and the set 
$$\displaystyle  S_{c} \ = \ \left\{ x \in \bar{\Omega} \mbox{  : } 1 \geq d_{c}
(x) \geq 1 - \mathcal{O}(c) \right\}$$ 
approximates the crack.  More precisely,  if $\displaystyle (u_{c},d_{c})$ is a minimizer of the Ambrosio-Tortorelli functional (\ref{atf}) then as $c \rightarrow 0$ the displacement $\displaystyle u_{c}$ converges (in some norm) to a displacement $u$, the set $\displaystyle S_{c}$ shrinks to a surface $S$ and $(u,S)$ is a minimizer of the Mumford-Shah functional 
 (\ref{ms2}).

The variable $d$  plays the role of a brittle damage variable, because it takes values in $[0,1]$ and also because 
it is coupled with the antiplane displacement $u$ through the term 
$$ \int_{\Omega} \left\{ \phi(d) \, \frac{1}{2} \, K \, \mid \nabla u \mid^{2} \right\} \,  $$
which represents the elastic energy of the body with elasticity coefficient 
$\displaystyle  \phi(d) K$. The function $\phi$ is taken as  a decreasing function from $[0,1]$ to $[0,1]$, such that $\phi(1) = 0$,  $\phi(0) = 1$. 

Focardi \cite{focardi} proved that there is  a Ambrosio-Tortorelli functional suitable for approximating 
the 3D Mumford-Shah functional (\ref{ms3}), namely: 
\begin{equation}
E_{c}(u,d) \ = \  \int_{\Omega} \left\{ \phi(d) \, w(\nabla \eu) \, + \,
\frac{1}{2} \gamma \, c \, 
  \mid \nabla d \mid^{2} \, + \, \frac{\gamma}{2 c} \, d^{2} \right\} \,  
  \label{atf3d}
\end{equation}
under certain growth conditions on the elastic energy function $w$.

\subsection{Quasistatic model, using Mielke's theory}

In this subsection we obtain an interpretation of a mathematical result of 
Giacomini \cite{giacomini},  which shows that models of damage based on 
the Ambrosio-Tortorelli functional have the important property of being
compatible with brittle damage from the energetic point of view. This is a
desirable feature of a model of brittle damage, as there are many "classical" 
models of brittle damage which allow the creation of a brittle crack (seen 
a concentrated total damaged region) with zero consumed  energy.

We shall look at  the equations coming from the force balance equation of 
Mielke (\ref{fmielke}) and  the Ambrosio-Tortorelli functional taken as the 
potential energy.  The state of the system is described by a pair $(\eu, d)$, 
where $\eu$ is the displacement and $d$ a scalar damage variable taking values 
in $[0,1]$.

 We shall take a dissipation metric  which is almost the same 
 as in Mielke and Roub\'{\i}\v{c}ek 
 model \cite{MR06b}, relation (2.5) (see also the discussion at the end of the
 section 2.2), which gives the dissipation functional 
 $$\mathcal{R}(\eu, d, \dot{\eu}, \dot{d}) \, = \, \int_{\Omega} \left\{\beta  \dot{d} \, + \, \chi_{1}(d) \, + \chi_{2}(\dot{d}) \right\}  $$
 The functions $\displaystyle \chi_{1}, \chi_{2}$ are indicator functions of convex sets:
 $$\chi_{1}(d) \, = \, \left\{ 
 \begin{array}{ccl}
 0 & , & \mbox{ if } d \in [0,1] \\
 +\infty & , & \mbox{ else }
 \end{array} \right. \quad , \quad 
 \chi_{1}(\dot{d}) \, = \, \left\{ 
 \begin{array}{ccl}
 0 & , & \mbox{ if } \dot{d} \in [0,+\infty) \\
 +\infty & , & \mbox{ else }
 \end{array} \right.
$$

Formally integrating by parts the force balance equation of Mielke 
(\ref{fmielke}), we arrive  to the evolution equations: 
\begin{equation}
 \left\{ 
 \begin{array}{rcl}
0 & = &   \frac{\partial }{\partial \, x_{i}} \, \left( \phi(d) \, K \,
\frac{\partial \, u }{\partial x_{i}} \right) \\ 
0 & \in &  \beta \,  - \,  \gamma \, c \, \Delta d  \, + \, \phi'(d ) \, \frac{1}{2} K \, \mid \nabla u  \mid^{2} \, + \, \frac{\gamma}{2 c} \, d   \, + \, \partial \, \chi_{2} (\dot{d })
\end{array}
\right. \quad  .
\label{miat}
\end{equation}
We  add the constraints $d \in [0,1]$, boundary and initial conditions.  The term 
$$- \, \phi'(d ) \, \frac{1}{2} K \, \mid \nabla u  \mid^{2} $$
is greater or equal than $0$, due to the fact that $\phi$ is decreasing, thus 
$\displaystyle \phi' \leq 0$.  This term represents the variation of the 
elastic energy density due to damage.

The paper \cite{giacomini} can be seen as an investigation o f
 the limit to the fracture model of the bulk damage model of Mielke and 
 Roub\'{\i}\v{c}ek, that is in the limit  when the damage variable equals 
 $0$ almost everywhere (therefore the value of the parameter $\beta$ is not
 important in the sense that  $a>0$ makes the same effect as 
 $\beta =0$).  
This result can be described as follows:  for any parameter $c$ let 
$\displaystyle q_{c} = (\eu_{c}, d_{c})$ denote an energetic solution 
associated with the Ambrosio-Tortorelli energy $\displaystyle E_{c}$ and 
dissipation distance $D$ coming from the dissipation metric $\mathcal{R}$. 
Then as $c$ converges to $0$, the evolution $\displaystyle q_{c}$ 
converges to an evolution $(\eu, S)$ of the energetic formulation of brittle 
fracture of Francfort, Marigo \cite{FMa} or Buliga \cite{bu3}.

From the point of view of mechanics fracture is a manifestation of concentrated 
damage. Therefore a good (bulk) damage model should have the property that 
it is not possible to produce arbitrarily concentrated damage with arbitrarily
small expense of energy. Such models are said to be compatible with brittle
fracture from the viewpoint of energy balance. 
 There are many models of brittle damage in use, not all of them 
compatible with brittle damage. The
mathematical result of Giacomini means that  the Ambrosio-Tortorelli functional  leads  to 
brittle damage models which are compatible with brittle fracture from the 
point of view of energy balance.

\subsection{Hamiltonian brittle damage}

We shall apply the generalized Hamiltonian approach to a functional of the 
Ambrosio-Tortorelli type. 

We take as state $Q = (\eu, d)$ the pair formed by the  displacement $\eu$ and 
the scalar  damage 
variable $d \in [0,1]$. The space of this pairs corresponds to the space $X$ 
from the general model.  

The dual variable, in the sense of Hamiltonian 
mechanics, is $P = (\ep, y) \in Y $, where $\ep$ is the momentum  and $y$ is a 
scalar  variable dual to $d$ (which will turn out to be 
linearly dependent on  $\displaystyle \dot{d}$). 

The space of all pairs $(Q,P)$ is a product of two symplectic vector spaces 
$\displaystyle \mathcal{B} = \mathcal{B}_{1} \times \mathcal{B}_{2}$. 
The space of non-dissipative variables 
$\displaystyle B_{1} \times B_{1}^{*}$ is a  space of pairs of 
(weak) functions $(\eu, \ep)$ defined over $\Omega \subset \mathbb{R}^{3}$. 
Therefore 
$\displaystyle \eu \in  B_{1}$ and $\displaystyle \ep \in B_{1}^{*}$, where 
$\displaystyle   B_{1}$ is a Banach space (for example a well chosen Sobolev space 
of functions over $\Omega$) and $\displaystyle B_{1}^{*}$ is its dual. 
The duality product is 
$$\langle \ep, \eu \rangle_{1} \, = \, \int_{\Omega} \ep \cdot \eu   $$
Similarly, the space of dissipative variables $(d, y)$ is 
 $\displaystyle \mathcal{B}_{2} \, = \, B_{2} \times
B_{2}^{*}$,  a  space of pairs of (weak) functions $(d, y)$ defined over 
$\Omega \subset \mathbb{R}^{3}$, with $\displaystyle   B_{2}$  another 
Banach space  of functions over $\Omega$,   $\displaystyle  B_{2}^{*}$ is its dual. The duality product is 
$$\langle y, d \rangle_{2} \, = \, \int_{\Omega} y \, d   $$

Let us define the  the Hamiltonian as: 
\begin{equation}
H(t, \eu, \ep, d, y) \, = \, \mathcal{E}(\eu,  d) \, + \, T(\ep, y) \, - \, \langle l(t) , \eu \rangle
\label{hengen}
 \end{equation}
where $\Psi$ is the free energy, $T$ the kinetic energy and $l(t)$ the external forces, seen as:
$$\langle l(t) , \eu \rangle \, = \, \int_{\Omega} \ef(t) \cdot \eu   \, + \, \int_{\Gamma} \bar{\ef}(t)  \cdot \eu$$
Here $\Gamma \subset \partial \Omega$ is that part of the boundary where surface forces $\bar{\ef}(t)$ are imposed at the moment $t$. 
The stored energy $\mathcal{E}$ is therefore: 
$$\mathcal{E}(t, \eu, d) \, = \, \Psi(\eu, d) \, - \langle l(t), \eu \rangle
\quad . $$
Displacements may be imposed on another part $\Gamma'$ of the boundary $\partial \Omega$. This is 
done by imposing that at every moment $t \in [0,T]$ the displacement $\eu(t)$ belongs to a 
subspace $\displaystyle B_{1}(t) \subset B_{1}$ of {\bf kinematically admissible} displacements.  

The expression of the free energy is the following:
\begin{equation}
\Psi(\eu,  d) \, = \, \int_{\Omega} \left[ \phi(d) w(\nabla \eu) \, + \,
\frac{1}{2} K \| \nabla d \|^{2} \, + \, \frac{1}{2} L \mid d \mid^{2} \right]   
\label{freee}
\end{equation}
which has a form analoguous with the one proposed by Stumpf and Hackl
\cite{sthk} formula (3.34). Here $\phi$ is a smooth, decreasing  function with values in the  interval $[0,1]$. 

The  kinetic energy has the form: 
\begin{equation}
T(\ep,  y) \, = \, \int_{\Omega} \left[ \frac{1}{2} b \mid y \mid^{2}  \, + \, \frac{1}{2 \rho}  \| \ep \|^{2}  \right]   \quad \quad .
\label{kine}
\end{equation}
The second term in the expression of the kinetic energy is just the usual kinetic energy expresses as a
function of momentum $\ep$, as it is usual in the Hamiltonian formalism. 
Similarly, $y$ is a momentum variable corresponding to $d$ and  $b$ is the 
scalar version of a  microinertia tensor (we use the same 
name as Stumpf and Hackl \cite{sthk} concerning the kinetic energy described in their formula (2.4)). 
We suppose that the constants  $K$, $L$ and $b$ are  positive. 

The dissipation function  is the same as in the previous section: 
\begin{equation}
\mathcal{R}(d, y, \dot{d}, \dot{y}) \, = \, \int_{\Omega} \left[ \chi_{[0,1]} (d) \, + 
\, \chi_{[0,+\infty)}(\dot{d}) \,
+ \, \beta \mid \dot{d} \mid \right]    \quad \quad .
\label{nonsvect}
\end{equation}

We shall find the equations satisfied by any  curve of evolution 
$(\eu,  \ep, d,  y) : [0,T] \rightarrow \mathcal{B}$ which is a solution of
the generalized Hamiltonian equations (\ref{mainsymp}), for the  
Hamiltonian (\ref{hengen}) and dissipation (\ref{nonsvect}). By using the expressions of the free energy (\ref{freee}) and kinetic energy (\ref{kine}), we obtain: 
\begin{equation}
 \langle \dot{\ep} , \hat{\eu} \rangle_{1} \, + \,  \langle D_{\eu} \Psi(\eu, d) , \hat{\eu} \rangle_{1} \, = \, 
\langle l(t) , \hat{\eu} \rangle_{1} \quad  \quad \forall \hat{\eu} \in B_{1}(t)
\quad ,
\label{h1}
\end{equation}

\begin{equation}
\langle \hat{\ep} , \dot{\eu} \rangle_{1} \,  - \,  \langle \hat{\ep} , D_{\ep} K(\ep, y)  \rangle_{1} \, = \, 0 
 \quad  \quad \forall \hat{\ep} \in B_{1}^{*}   \quad  .
\label{h2}
\end{equation}
There are two more equations, for the evolution of $d$ and $y$. 
Due to the non smooth dissipation, these are in fact expressed as 
subdifferential inequalities: for almost any $t \in [0,T]$ 
$\displaystyle d(t,x) \in [0,1]$ and $y(t,x) \in [0, +\infty)$ for 
almost every $x \in \Omega$, at any $t \in [0,T]$ 
the displacement $\eu(t)$ is kinematically admissible, i.e. $\displaystyle \eu(t) \in B_{1}(t)$,  and moreover 
 for any $\displaystyle  \hat{d} \in B_{2}$, such that $\displaystyle
\hat{d}(x) + \dot{d} \geq 0$ for almost every $x \in \Omega$, and for any 
$\displaystyle \hat{y} \in B_{2}^{*}$  we have: 
\begin{equation}
\beta \, \int_{\Omega} \left[ \mid \dot{d} + \hat{d} \mid \, - \, \mid \dot{d}
\mid \right]   \, \geq \,
\label{h3}
\end{equation} 
$$\geq \, \int_{\Omega} \left[ \hat{y} (\dot{d} - b y) \, - \, 
(L d + \phi'(d) w(\nabla\eu) + \dot{y}) \hat{d} \, - \, 
K \nabla d \, \nabla \hat{d} \right]      \quad \quad . $$

The equation (\ref{h1}) gives the usual momentum balance: for any $\hat{\eu}$ kinematically admissible 
we have 
$$\int_{\Omega} \left[ - \, \dot{p} \cdot \hat{\eu} \, - \phi(d) \, Dw(\nabla \eu) : \nabla \hat{\eu} \right]   \, = \, 
\int_{\Omega} \ef(t) \cdot \hat{\eu}   \, + \, \int_{\Gamma} \bar{\ef}(t) \cdot \hat{\eu}$$
Denote by $\displaystyle \eS = Dw(\nabla \eu)$ the stress variable given by the elastic energy $w$. Integration by parts leads us to a balance equation and boundary conditions: 
\begin{equation}
div \, \left( \phi(d) \, \eS \right) \, + \, \ef(t) \, = \, \dot{\ep} \quad \mbox{ in }Ê\Omega
\label{h1loceq}
\end{equation}
\begin{equation}
\phi(d) \, \eS \en \, = \, \bar{\ef}(t)   \mbox{ on } \Gamma \,  , \, \phi(d) \, \eS \en \, = \, 0 
   \mbox{ on } \partial \Omega \setminus \left( \Gamma \cup \Gamma'\right) \, , \, \eu \, = \, \eu_{0}(t)  
   \, \mbox{ on } \Gamma'   \quad \quad .
\label{h1locb}
\end{equation}
Equation (\ref{h2}) gives us the momentum $\ep$ as function of 
$\displaystyle \dot{\eu}$: 
\begin{equation}
\ep \, = \, \rho \, \dot{\eu}   \quad \quad . 
\label{h2loc}
\end{equation}
Equation (\ref{h3}) is equivalent to the following two relations:  
\begin{equation}
\dot{d} \, = \, b y
\label{nq}
\end{equation}
and for all $\displaystyle  \hat{d} \in B_{2}$, such that $\displaystyle
\hat{d}(x) + \dot{d} \geq 0$ for almost every $x \in \Omega$
\begin{equation}
\beta \, \int_{\Omega} \left[ \mid \dot{d} + \hat{d} \mid \, - \, \mid \dot{d}
\mid \right]   \, \geq \,
\label{newq}
\end{equation}
$$ \geq \, - \int_{\Omega}  \left[  
(L d + \phi'(d) w(\nabla\eu) + \dot{y}) \hat{d} \, + \, 
K \nabla d \, \nabla \hat{d} \right]      \quad \quad . $$
Let $\displaystyle S: [0, +\infty) \rightarrow 2^{\mathbb{R}}$ be the
multivalued function defined by: 
$$S(v) = \left\{ 
\begin{array}{lcl}
\beta & , & v > 0 \\
(-\infty, \beta] & , &  v =  0
\end{array} \right. \quad \quad . $$
The function $S$ is the subdifferential of a convex function. By using the definition
of $S$ and relation (\ref{nq})  we obtain the following equivalent form 
of the  inequality (\ref{newq}): for almost every $x \in \Omega$ we have: 
\begin{equation}
- \, \left( \dot{y} \, + \, L d \, + \, \phi'(d) w(\nabla\eu) \, - \, K \Delta d
\right) \, \in \, S(y)  \quad \quad . 
\label{newqq}
\end{equation}
We may add the boundary condition (which is not strictly speaking a consequence
of the formalism): on $\partial \Omega$ we have: 
\begin{equation}
y \, \geq \, 0 \, , \quad - K \frac{d}{d \en} d \in S(y)  \quad \quad .
\end{equation}

In the particular case of a functional of the Ambrosio-Tortorelli type  (\ref{atf3d}) we may take:  
$$ K \, = \, \gamma \, c \quad , \quad L \, = \, \frac{\gamma}{c} \quad , \quad 
b \, = \, \gamma \, c  $$
 The function $\phi$ which enters in the expression of 
the free energy is chosen as in the Ambrosio-Tortorelli functional.  With this 
choice of constants we obtain from (\ref{newqq}) and (\ref{nq}) the differential
inclusion: 
$$- \, \left( \ddot{d} \, + \, \gamma^{2} d \, + \, 
\gamma c \, \phi'(d) w(\nabla\eu) \, - \, \gamma^{2} c^{2} \Delta d
\right) \, \in \, \gamma c \, S(\dot{d})  \quad \quad . $$
This inclusion suggests that in this model there is a maximal speed of
propagation of damage of order $\displaystyle \gamma \sqrt{1 + c^{2}}$.

\end{document}